\documentclass[a4paper]{article}

\usepackage{amsmath, amsfonts, amsthm, amssymb}
\usepackage[T1]{fontenc}
\usepackage{ae}

\newcommand\HH{\mathcal{H}}
\newcommand\CC{\mathbb{C}}
\renewcommand\S{{S}}
\newcommand\f{{f}}
\newcommand\gone{{g_1}}
\newcommand\gtwo{{g_2}}

\newcommand\QQ{\mathbb{Q}}
\newcommand\N{N}
\newcommand\NN{\mathcal{N}}
\newcommand\RR{\mathbb{R}}
\newcommand\R{\mathcal{R}}
\newcommand\TT{\mathcal{T}}
\newcommand\ZZ{\mathbb{Z}}
\newcommand\ZZp{\ZZ_{>0}}
\newcommand\PP{\mathbb{P}}
\newcommand\xx{\mathbf{x}}
\newcommand\dd{\,\mathrm{d}}
\newcommand\Pthree{{\PP^3}}
\newcommand{\base}[7]{\xi^{({#1},{#2},{#3},{#4},{#5},{#6},{#7})}}
\newcommand{\congr}[3]{{#1} \equiv {#2}\ (\mathrm{mod}\ {#3})}
\newcommand\kell{{k_\ell}}
\newcommand\tell{\tau_\ell}
\newcommand\tone{\tau_1}
\newcommand\ttwo{\tau_2}
\newcommand\ep{\epsilon}
\newcommand\lle{{\ll_\ep}}
\newcommand\re{{\Re e}}
\newcommand\NU{{N_{U,H}}}
\newcommand\tpsi{{\widetilde \psi}}
\newcommand\Res{\mathrm{Res}}
\newcommand{\Esix}{{\mathbf E}_6}
\newcommand{\Dfour}{{\mathbf D}_4}
\newcommand{\Dfive}{{\mathbf D}_5}
\newcommand{\Fell}{\xi_\ell^3\xi_4^2\xi_5}
\newcommand{\Ftwo}{\xi_2}
\newcommand{\Fone}{\xi_1^2\xi_3}
\newcommand{\Tone}{\tone^3\Fone}
\newcommand{\Ttwo}{\ttwo^2\Ftwo}
\newcommand{\Tell}{\tell\Fell}
\newcommand{\xell}{{\xi_\ell}}
\newcommand{\Cone}{\xi_2\xi_3\xi_\ell\xi_4\xi_5\xi_6}
\newcommand{\Ctwo}{\xi_1\xi_3}
\newcommand{\Cell}{\xi_4\xi_5\xi_6}
\newcommand{\hcfone}{\hcf(\tone,\Cone)=1}
\newcommand{\hcftwo}{\hcf(\ttwo,\Ctwo)=1}
\newcommand{\hcfell}{\hcf(\tell,\Cell)=1}
\newcommand{\tS}{{\widetilde \S}}
\newcommand{\anti}{-K_\tS}
\newcommand{\hcfrho}{{\hcf(\varrho,\kell\Fell) = 1}}
\newcommand{\congrrho}{{\congr{\ttwo}{\varrho\tone\xi_1}{\kell\Fell}}}
\newcommand{\congrrhotwo}{{\congr{-\varrho^2\xi_2}{\tone\xi_3}{\kell\Fell}}}
\newcommand{\boundrho}{{0 < \varrho \le \kell\Fell}}
\renewcommand{\le}{\leqslant}
\renewcommand{\ge}{\geqslant}

\DeclareMathOperator{\Pic}{Pic}
\DeclareMathOperator{\Spec}{Spec}
\DeclareMathOperator{\hcf}{gcd}
\DeclareMathOperator{\Vol}{Vol}
\DeclareMathOperator{\rk}{rk}

\newtheorem{theorem}{Theorem}
\newtheorem{lemma}[theorem]{Lemma}
\newtheorem{proposition}[theorem]{Proposition}
\theoremstyle{definition}

\numberwithin{equation}{section}

\begin{document}

\title{Manin's conjecture\\for a certain singular cubic surface}
\author{Ulrich Derenthal\thanks{Mathematisches Institut der
    Universit\"at G\"ottingen, Bunsenstr. 3--5, D-37073 G\"ottingen,
    Germany. email: derentha@math.uni-goettingen.de}} \maketitle

\begin{abstract}
  This paper contains a proof of Manin's conjecture for the singular
  cubic surface $\S \subset \Pthree$ with a singularity of type
  $\Esix$, defined by the equation $x_1x_2^2+x_2x_0^2+x_3^3=0$. If $U$
  is the open subset of $\S$ obtained by deleting the unique line from
  $\S$, then the number of rational points in $U$ with anticanonical
  height bounded by $B$ behaves asymptotically as $cB(\log B)^6$,
  where the constant $c$ agrees with the one conjectured by Peyre.
\end{abstract}

\section{Introduction}

Let $\S$ be a Del Pezzo surface, possibly singular. If its set of
rational points $\S(\QQ)$ is non-empty, then it is dense in the Zariski
topology, and it is natural to ask for the asymptotic behaviour of the
number of rational points of bounded height.

If $\S$ contains lines defined over $\QQ$ then the number of rational points
of bounded height on these dominates the number of rational points 
in their complement $U\subset \S$.
Therefore, one is mainly interested in rational points in $U$.

In this situation, Manin has formulated a
far-reaching conjecture for the behaviour of the counting function
\[
\NU(B):=\#\{x \in U(\QQ) \mid H(x) \le B\},
\]
where $H$ is an exponential height \cite{MR89m:11060}, \cite{MR1032922}.
In the special case, when $H$ is associated to the
anticanonical embedding of $\S$, the conjecture states that 
\[
\NU(B) \sim c B (\log B)^r,
\] 
where $r$ depends only on the geometry of $\S$. 
A conjectural interpretation of $c$ was
later given by Peyre \cite{MR1340296}.

In recent years, Manin's conjecture has been proved for several classes of
algebraic varieties. Batyrev and Tschinkel proved it for toric varieties
\cite{MR1620682}, in particular for all Del Pezzo surfaces of degree $\ge 6$.
For smooth Del Pezzo surfaces of degree 5, it was proved by de la Bret{\`e}che
\cite{MR2003m:14033} in the split case (i.e., with geometric Picard group
defined over $\QQ$) and later together with Fouvry in several non-split cases
\cite{MR2099200}. In degree 4, it has been established recently in the case of
a split surface with a $\Dfive$-singularity \cite{math.NT/0412086} and a
non-split surface with a $\Dfour$-singularity \cite{math.NT/0502510}. In
degree 3, apart from the singular toric surfaces, Manin's conjecture and its
refinement by Peyre have not yet been established. The best progress so far
are upper and lower bounds of the right order of magnitude proved by
Heath-Brown in case of the Cayley cubic \cite{0981.32025}, and by Browning for
a surface with a singularity of type $\Dfour$ \cite{math.NT/0404245}.

Manin's conjecture seems to be more difficult for
surfaces of lower degree or with fewer singularities.
The best general upper bound for $\NU(B)$ for a smooth cubic surface is 
$B^{4/3+\ep}$, given by Heath-Brown \cite{MR98h:11083}, \cite{MR2000a:11042}.

\bigskip

We prove Manin's conjecture in the case of a cubic surface
with a singularity of type $\Esix$:
\begin{equation}\label{eq:surface}
  \S = \{\xx = (x_0:x_1:x_2:x_3) \in \Pthree \mid \f(\xx) =
x_1x_2^2+x_2x_0^2+x_3^3 = 0\}
\end{equation}
Here, $\Esix$ refers to the Dynkin diagram
describing the intersections of exceptional divisors 
in the minimal desingularization $\tS$ of $\S$.

\begin{theorem}\label{thm:main}
Let $\S$ be the cubic surface with an $\Esix$-singularity as above and let $H$
be the associated height. Let $U := \S \setminus \ell$, where $\ell =
\{x_2=x_3 = 0\}$ is the unique line on $\S$.
  
Then 
\[\NU(B) = c_{\S,H} B Q(\log B) + O(B(\log B)^2),\] 
where $Q$
is a monic polynomial of degree $6$, and the leading constant
$c_{\S,H}$ is the one conjectured by Peyre \cite{MR1340296}.
\end{theorem}

The invariants appearing in Manin's conjecture and Peyre's constant $c_{\S,H}$
are calculated in Section~\ref{sec:manins_conjecture}.

\

The proof follows the strategy of de la Bret{\`e}che and Browning
\cite{math.NT/0412086} and uses the universal torsor.  The passage to the
universal torsor was introduced by Colliot-Th\'el\`ene and Sansuc in their
investigations of the Hasse principle on Del Pezzo surfaces \cite{MR54:2657},
\cite{MR89f:11082}. Salberger proposed using torsors in the study of rational
points of bounded height \cite{MR1679841}, see also \cite{MR1679842},
\cite{MR2029862}.  For our particular surface, the universal torsor has been
calculated by Hassett and Tschinkel using the Cox ring; the torsor is a
hypersurface in the $10$-dimensional affine space \cite{MR2029868}.  In
addition to the equation defining the torsor we need to derive certain
coprimality conditions between the coordinates. In
Section~\ref{sec:universal_torsor} we compute the torsor and determine these
conditions following the more direct approach of Heath-Brown, Browning and de
la Bret{\`e}che \cite{0981.32025}, \cite{math.NT/0404245},
\cite{math.NT/0412086}.

The next step is to count the number of integral points on the
universal torsor subject to certain bounds, given by lifting the height
function to the torsor, and satisfying the coprimality conditions.
For three of the ten variables on the torsor, this summation is done
in Section~\ref{sec:summations} by elementary methods from analytic
number theory. The summation over the last seven variables, 
completing the proof of Theorem~\ref{thm:main},
is carried out in Section~\ref{sec:main_thm_proof}.

\

During the final preparation of this paper, 
I was informed by M. Joyce 
(Brown University) that similar results towards 
Manin's conjecture for certain cubic surfaces 
will appear in his thesis.

\

\paragraph{Acknowledgments.}
Most of this paper was written while I was visiting B. Hassett at
Rice University in February 2005. I am very grateful to 
him for the invitation and stimulating discussions.
I thank T. D. Browning for helpfully explaining the details of his paper
\cite{math.NT/0412086}. I thank S. Elsenhans and E. Kappos for many useful 
suggestions. 
This research has been partially supported by the 
Studienstiftung des Deutschen Volkes.

\section{Manin's Conjecture}\label{sec:manins_conjecture}

In this section, we calculate the invariants
appearing in Manin's conjecture and its refinement by Peyre. 
We use the notation from \cite{MR1340296}.

\begin{lemma}\label{lem:manin_conjecture}
  Manin's conjecture predicts in case of $\S$ as defined in
  \eqref{eq:surface}:
  \[\NU(B) \sim c_{\S,H}B(\log B)^6,\] where $c_{\S,H} =
  \alpha(\tS)\beta(\tS)\omega_H(\tS)$ with
  \begin{equation*}
    \begin{split}
      \alpha(\tS) &:= \frac{1}{6! \prod_i \lambda_i} = \frac{1}{6!\cdot 2
        \cdot 3 \cdot 4 \cdot 3 \cdot 4 \cdot
        5 \cdot 6}=\frac 1{6220800},\\
      \beta(\tS) &:= 1,\\
      \omega_H(\tS) &:= \omega_\infty \cdot \omega_0,
    \end{split}
  \end{equation*}
  and
  \begin{equation*}
    \begin{split}
      \lambda &= (\lambda_1, \lambda_2, \lambda_3, \lambda_\ell,
      \lambda_4,
      \lambda_5, \lambda_6) := (2,3,4,3,4,5,6),\\
      \omega_\infty &:= 6 \int\int\int_{\{(u,v,w)\in
        \RR^3 \mid |tv^3| \le 1,
        |t^2+u^3|\le 1, 0 \le v \le 1, |uv^4| \le 1\}} 1 \dd t \dd u \dd v,\\
      \omega_0 &:= \prod_p\omega_p,\\
      \omega_p &:= \left(1-\frac 1 p\right)^7\left(1+\frac 7 p + \frac
        1 {p^2}\right).
    \end{split}
  \end{equation*}
\end{lemma}
\begin{proof}
The line  $\ell$ and all exceptional divisors $F_1, \dots, F_6$
of the desingularization $\tS \rightarrow \S$ are defined over $\QQ$.
Thus, $\rk \Pic(\tS) = 7$. 
By \cite{MR2029868}, the effective cone is generated by the classes of $\ell$
and the exceptional divisors. 
The dual cone of nef divisors is also simplicial (i.e., the number of
generators is equal to $\rk \Pic(\tS)$), and 
\[
\anti = 2F_1 + 3F_2 + 4F_3+3\ell+4F_4+5F_5+6F_6.
\] 
Note that $\lambda=\anti\in \Pic(\tS) \cong \ZZ^7$, in the basis 
$\{F_1, F_2, F_3,\ell, F_4, F_5, F_6\}$.

  Using the definition, we calculate
  \[\alpha(\tS) = \Vol\bigg\{(t_1,t_2,t_3,t_\ell,t_4,t_5,t_6) \in
      \RR^7_{\ge 0} \bigg| \sum_{i \in \{1,2,3,\ell,4,5,6\}} \lambda_i
      t_i = 1\bigg\} =\frac 1{6!\cdot \prod_i \lambda_i}.\]
The surface $\S$ is split over $\QQ$, so that
\[
\beta(\tS) = \#H^1(\QQ,\Pic(\tS))= 1.
\]
By definition, 
\[
\omega_H(\tS) = \lim_{s \to 1}((s-1)^{\rk
    \Pic(\tS)} L(s,\Pic(\tS))) \omega_\infty \prod_p \frac
  {\omega_p}{L_p(1,\Pic(\tS))},
\] where, in our case,
  \[\lim_{s \to 1}((s-1)^{\rk \Pic(\tS)} L(s, \Pic(\tS)) = \lim_{s \to
    1}((s-1)^7\zeta(s)^7) = 1\]
  and \[L_p(1, \Pic(\tS))^{-1} = (1 - 1/p)^7.\]
  
  We use Peyre's method \cite{MR1340296} to compute $\omega_\infty$
  and parametrize the points by writing $x_1$ as a function of
  $x_0,x_2,x_3$. Since $\xx = -\xx$ in $\Pthree$, we may assume $x_2
  \ge 0$. Since $\frac{\dd}{\dd x_1}\f = x_2^2$, the Leray form
  $\omega_L(\tS)$ is given by $x_2^{-2} \dd x_0 \dd x_3 \dd x_2$,
  and we obtain $\omega_\infty(\tS)$ from
  \[\int\int\int_{\{|x_0|\le 1,|x_2^{-2}(x_2x_0^2+x_3^3)|\le 1,0 \le x_2
    \le 1,|x_3|\le 1\}} x_2^{-2} \dd x_0 \dd x_3 \dd x_2,\]
  using the transformations \[x_0 = tx_2^{1/2}, \qquad x_3 =
  ux_2^{2/3}, \qquad x_2 = v^6.\]
The calculation of $\omega_p$ is done as in
Lemma 1 of \cite{math.NT/0412086} and we omit it here.
\end{proof}

\section{The Universal Torsor}\label{sec:universal_torsor}

The universal torsor $\TT =
\Spec(\QQ[\xi_1,\xi_2,\xi_3,\xi_\ell,\xi_4,\xi_5,\xi_6,\tau_1,\tau_2,\tau_\ell]
/ (T(\xi_i,\tau_i)))$ is given by the equation
\begin{equation}\label{eq:torsor}
  T(\xi_i,\tau_i) = \Tell+\Ttwo+\Tone = 0
\end{equation}
and the map $\Psi: \TT \to \S = \Spec(k[x_0, \dots, x_3]/(\f(\xx)))$
defined by
\begin{equation}\label{eq:substitutions}
  \begin{split}
    \Psi^*(x_0) &= \base 1 2 2 0 1 2 3 \tau_2\\
    \Psi^*(x_1) &= \tau_\ell\\
    \Psi^*(x_2) &= \base 2 3 4 3 4 5 6\\
    \Psi^*(x_3) &= \base 2 2 3 1 2 3 4 \tau_1
  \end{split}
\end{equation}
where we use the notation $\base {n_1}{n_2}{n_3}{n_\ell}{n_4}{n_5}{n_6} :=
\xi_1^{n_1}\xi_2^{n_2}\xi_3^{n_3} \xi_\ell^{n_\ell}
\xi_4^{n_4}\xi_5^{n_5}\xi_6^{n_6}$. Note that $\Psi^*(x_2) = \xi^\lambda$
with $\lambda \in \ZZ^7$ as in Lemma~\ref{lem:manin_conjecture}.

We want to establish a bijection between rational points on the surface $\S$
and integral points on the torsor $\TT$ which are subject to certain
coprimality conditions.

More precisely, the coprimality conditions can be summarized in the following
table, where a ``$-$'' means that the two variables are coprime, and a
``$\times$'' that they may have common factors. For a variable combined with
itself, ``$-$'' means that each prime occurs at most once, and
``$\times$'' means that it may occur more often.
\begin{equation*}
\begin{array}{c|ccccccc|ccc}
& \xi_1 & \xi_2 & \xi_3 & \xi_\ell & \xi_4 & \xi_5 & \xi_6 & \tau_1 & \tau_2 &
\tau_\ell\\
\hline
\xi_1 &\times&-&\times&-&-&-&\times&\times&-&-\\
\xi_2 &-&-&-&-&-&-&\times&-&\times&-\\
\xi_3 &\times&-&-&-&-&-&\times&-&-&-\\
\xi_\ell &-&-&-&\times&\times&\times&\times&-&-&\times\\
\xi_4 &-&-&-&\times&-&-&\times&-&-&-\\
\xi_5 &-&-&-&\times&-&-&\times&-&-&-\\
\xi_6 &\times&\times&\times&\times&\times&\times&\times&-&\times&-\\
\hline
\tau_1 &\times&-&-&-&-&-&-&\times&\times&\times\\
\tau_2 &-&\times&-&-&-&-&\times&\times&\times&\times\\
\tau_\ell &-&-&-&\times&-&-&-&\times&\times&\times\\
\end{array}
\end{equation*}

Later, we will refer to the
\begin{equation}\label{eq:hcf_xi}
  \text{coprimality conditions between $\xi_1, \dots, \xi_6$ 
    as given in the table}
\end{equation}
Because of the torsor equation $T$, we can write the coprimality
conditions for $\tau_i$ equivalently as
\begin{equation}\label{eq:hcf_tone}
  \hcfone
\end{equation}
and
\begin{equation}\label{eq:hcf_ttwo_tell}
  \hcftwo,\qquad \hcfell.
\end{equation}

The goal of this section is the following result:

\begin{proposition}\label{prop:bijection}
  The map $\Psi$ induces a bijection between \[\TT_1 := \{(\xi_i,\tau_i) \in
  \TT(\ZZ) \mid \text{\eqref{eq:hcf_xi}, \eqref{eq:hcf_tone},
    \eqref{eq:hcf_ttwo_tell} hold}, \xi_i > 0\}\] and $U(\QQ) \subset
  \S(\QQ)$.
\end{proposition}

The proof of this is split into two parts. First, we establish a similar
bijection with slightly different coprimality conditions:

\begin{lemma}\label{lem:bijection_different_coprime}
  Let $\TT_2$ be set of all $(\xi_i,\tau_i) \in \TT(\ZZ)$ such that
  $\xi_i > 0$ and the coprimality conditions described by the table
  hold, except that the conditions
  \begin{equation}\label{eq:new_coprim}
    \hcf(\xi_3,\tone)=1\text{ and }\hcf(\xi_6,\tone)=1
  \end{equation}
  in the table are replaced by 
  \begin{equation}\label{eq:old_coprim}
    |\mu(\xi_1)|=1\text{ and }\hcf(\xi_1,\xi_3)=1.
  \end{equation}
  Then the map $\Psi$ induces a bijection between $\TT_2$ and $U(\QQ) \subset
  \S(\QQ)$.
\end{lemma}

\begin{proof}
  Using the method of \cite{math.NT/0412086}, we now show
  that the coprimality conditions lead to a bijection.
  
  We go through a series of coprimality considerations and replace the
  original variables by products of new ones which fulfill certain conditions.
  When doing this, the new variables will be uniquely determined.
  
  Since $\xx = -\xx$, and $x_2 = 0$ is equivalent to $\xx \in \ell$, we can
  write each $\xx \in U(\QQ)$ uniquely such that $x_i \in \ZZ$, $x_2 > 0$, and
  $\hcf(x_i) = 1$.
\begin{itemize}
\item Note that $x_2|x_3^3$. Write $x_2=y_1y_2^2y_3^3$ with $y_i \in \ZZp$,
  where each triple occurrence of a prime factor of $x_3$ is put in $y_3$ and
  each double occurrence in $y_2$, so that $y_1,y_2,y_3$ are unique if we
  assume $|\mu(y_1y_2)|=1$.  Then $x_3 = y_1y_2y_3z$ must hold for a suitable
  $z \in \ZZ$. Substituting into $\f$ and dividing by $y_1y_2^2y_3^3$ gives
  \[\f_1(x_0,x_1,y_1,y_2,y_3,z) = x_1y_1y_2^2y_3^3 + x_0^2 + y_1^2y_2z^3 = 0.\]

\item Now $y_1y_2|x_0^2$, and since $|\mu(y_1y_2)|=1$, we have $y_1y_2|x_0$.
  Write $x_0=y_1y_2w$, where $w\in \ZZ$. Substituting and dividing by
  $y_1y_2$, we obtain 
  \[\f_2(x_1,y_1,y_2,y_3,z,w) = x_1y_2y_3^3+w^2y_1y_2+y_1z^3 = 0.\]

\item Since $y_2|y_1z^3$ and $|\mu(y_1y_2)|=1$, we must have $y_2|z$. Write
  $z=y_2z'$, where $z' \in \ZZ$, and obtain, after dividing by $y_2$, the
  relation
  \[\f_3(x_1,y_1,y_2,y_3,w,z') = x_1y_3^3+w^2y_1+y_1y_2^2z'^3 = 0.\]

\item Since $y_1$ divides our original variables $x_0,x_2,x_3$, it cannot
  divide $x_1$. Together with $|\mu(y_1)|=1$, the fact $y_1|x_1y_3^3$ implies
  $y_1|y_3$. Write $y_3 = y_1y_3'$, where $y_3' \in \ZZp$ and obtain
  \[\f_4(x_1,y_1,y_2,w,z',y_3') = x_1y_1^2y_3'^3+w^2+y_2^2z'^3 = 0\]
\item Let $a=\hcf(y_3',z') \in \ZZp$ and write $y_3'=ay_3''$, where
  $y_3''\in \ZZp$ and $z'=az''$, where $z'' \in \ZZ$. This gives
  \[\f_5(x_1,y_1,y_2,w,z'',y_3'',a) = x_1y_1^2y_3''^3a^3+w^2+y_2^2z'a^3= 0.\]
\item Now $a^3|w^2$. Writing $a=\xi_6^2\xi_2$, where $\xi_2,\xi_6 \in \ZZp$
  with $|\mu(\xi_2)|=1$, gives $w=w'\xi_6^3\xi_2^2$, where $w'\in \ZZ$ leading
  to the equation
  \[\f_6(x_1,y_1,y_2,z'',y_3'',w',\xi_2,\xi_6) = 
  x_1y_1^2y_3''^3+w'^2\xi_2+y_2^2z''^3=0.\]
  
\item Let $\xi_5 = \hcf(y_3'',w') \in \ZZp$ and write $y_3''=\xell\xi_5$,
  where $\xell \in \ZZp$ and $w'=w''\xi_5$, with $w'' \in \ZZ$.  Then
  \[\f_7(x_1,y_1,y_2,z'',w'',\xi_2,\xell,\xi_5,\xi_6) = 
  x_1y_1^2\xell^3\xi_5^3 + w''^2\xi_2\xi_5^2 + y_2^2z''^3 = 0.\]
  
\item Since $\hcf(y_3'',z'') = 1$, also $\hcf(\xell\xi_5,z'') = 1$.
  Therefore, $\xi_5^2|y_2^2$, which means $\xi_5|y_2$, and we write
  $y_2=\xi_1\xi_5$, with $\xi_1\in \ZZp$. We obtain 
  \[\f_8(x_1,y_1,z'',w'',\xi_1,\xi_2,\xell,\xi_5,\xi_6) =
  x_1y_1^2\xell^3\xi_5+w''^2\xi_2+\xi_1^2z''^3 = 0.\]
  
\item Let $\xi_3 = \hcf(w'',y_1) \in \ZZp$. Since
  $|\mu(y_1y_2)|=1$, we have $\hcf(\xi_1,\xi_3) = 1$. Therefore,
  $\xi_3|z''^3$ and even $\xi_3|z''$.  Write $w''=\ttwo \xi_3$, where
  $\ttwo \in \ZZ$, $y_1=\xi_4\xi_3$, where $\xi_4 \in \ZZp$, and
  $z''=\tone \xi_3$, where $\tone \in \ZZ$. Replacing $x_1=\tell$, where
  $\tell \in \ZZ$, we get
  \[\f_9(\xi_1,\xi_2,\xi_3,\xell,\xi_4,\xi_5,\xi_6,\tone,\ttwo,\tell)
  = \tell \xell^3\xi_4^2\xi_5+\ttwo^2\xi_2+\tone^3\xi_1^2\xi_3 = 0.\]
  This is the torsor equation $T(\xi_i,\tau_i)$ as in \eqref{eq:torsor}.
\end{itemize}

The substitutions lead to $x_0, \dots, x_3$ in terms of $\xi_i,\tau_i$ as in
\eqref{eq:substitutions}. Conversely, it is easy to check that each
$(\xi_i,\tau_i)$ satisfying $T$ is mapped by $\Psi$ to a point $\xx \in
\S(\QQ)$. Note that $\xi_i \in \ZZp$ and $\tau_i \in \ZZ$.  Furthermore, the
coprimality conditions we introduced impose the following conditions on
$\xi_i, \tau_i$:
\begin{equation*}
  \begin{split}
    &|\mu(y_1y_2)|=|\mu(\xi_1\xi_3\xi_4\xi_5)|=1,\qquad |\mu(\xi_2)|=1,
    \qquad \hcf(\ttwo,\xi_4)=1,\\
    &\hcf(y_3'',z'') = \hcf(\xell\xi_5,\tone\xi_3) = 1,\qquad
    \hcf(\xell,w'')
    = \hcf(\xell,\ttwo\xi_3) = 1.\\
  \end{split}
\end{equation*}
The condition $\hcf(x_i) = 1$ is equivalent to
$\hcf(\tell,\xi_1\xi_2\xi_3\xi_4\xi_5\xi_6)=1$.
  
  We obtain $\hcf(\xi_2,\xi_3)=1$ in the following manner: If
  $p|\xi_2,\xi_3$ for some prime $p$, then $p|\tell\Fell$ by the
  torsor equation $T$. On the other hand, a divisor of $\xi_3$ cannot
  divide any of the factors by the coprimality conditions we found.
  Similarly, we conclude $\hcf(\xi_3,\ttwo) = \hcf(\xi_1,\ttwo) =
  \hcf(\xi_2,\xi_5) = \hcf(\xi_5,\ttwo) = 1$.
  
  Finally, if a prime $p$ divides two of $\xi_2,\xi_4,\tone$, we see
  using $T$ that $p$ must divide all of them. Since $|\mu(\xi_2)|=1$ and
  $p^2|\Tell+\Tone$, we conclude $p|\ttwo$ which is impossible since
  $\hcf(\ttwo,\xi_4)=1$. Therefore, $\xi_2,\xi_4,\tone$ must be
  pairwise coprime. In the same way we derive that no two of
  $\xi_1,\xi_2,\xell$ have a common factor.
  
  It is easy to check that we cannot derive any other 
  coprimality condition because we could construct a solution to $T$
  violating it.
  
  Note that the conditions on $(\xi_i,\tau_i)$ are exactly the ones
  given in the definition of $\TT_2$. Since in every step the newly
  introduced variables are uniquely determined, we have established a
  bijection between $U(\QQ)$ and $\TT_2$.
\end{proof}

The following Lemma is the second step towards the proof of
Proposition~\ref{prop:bijection}:

\begin{lemma}\label{lem:bijection_different_correct}
  There is a bijection between $\TT_1$ and $\TT_2$.
\end{lemma}

\begin{proof}
  Given a point $(\xi'_i,\tau'_j) \in \TT_2$ violating
  \eqref{eq:new_coprim}, we could replace a common prime factor
  $p$ of $\xi'_3, \xi'_6$ and $\tau'_1$ by powers of $p$ as factors of
  $\xi'_1$ and possibly $\xi'_3$ such that \eqref{eq:new_coprim} holds.
  This way, we obtain a point $(\xi_i,\tau_j) \in \TT_1$. This should
  be done in a way such that $\Psi$ maps $(\xi_i,\tau_j)$ and
  $(\xi'_i,\tau'_j)$ to the same point $\xx \in U(\QQ)$, and such that
  we have an inverse map, taking care of the conditions
  \eqref{eq:old_coprim}.
  
  Let $(\xi_i, \tau_j) \in \TT_1$ and $(\xi'_i, \tau'_j) \in \TT_2$. Decompose
  the coordinates into their prime factors: Let 
  \[\xi_i = \prod_p p^{n_{ip}},\quad
  \tau_j = \prod_p p^{m_{jp}}\qquad\text{ and }\qquad \xi'_i = \prod_p
  p^{n'_{ip}},\quad \tau'_j = \prod_p p^{m'_{jp}},\] where $i \in
  \{1,2,3,\ell,4,5,6\}$ and $j \in \{1,2,\ell\}$.  Note that
  \eqref{eq:new_coprim} translates to $(n_{ip},m_{jp})$ fulfilling
  $n_{3p}=n_{6p}=0$ or $m_{1p}=0$, and that \eqref{eq:old_coprim} means that
  $(n'_{ip},m'_{jp})$ must fulfill $n'_{1p}+n'_{3p}\le 1$. Furthermore,
  $n_{3p}, n'_{3p} \in \{0,1\}$ always holds.
  
  Define the map 
\[\begin{array}{cccc}
\Phi':&  \TT_2         & \to     &  \TT_1 \\ 
      &(\xi'_i,\tau'_i)& \mapsto & (\xi_i,\tau_j),
\end{array}\] 
where $n_{ip} := n'_{ip}$ for $i \in \{2,\ell,4,5\}$ and
  $m_{jp} := m'_{jp}$ for $j \in \{2,\ell\}$, and the values of
  $n_{1p},n_{3p},n_{6p},m_{1p}$ depend on the size of $n'_{6p}$ compared to
  $m'_{1p}$, whether $m'_{1p}$ is even or odd, and whether $n'_{3p}$ is $0$ or
  $1$:
\begin{itemize}
\item 
  If $m'_{1p} = 2k+1, n'_{6p} \ge k+1, n'_{3p} = 0$, then
  \[(n_{1p},n_{3p},n_{6p},m_{1p}) :=
  (n'_{1p}+3k+1,1,n'_{6p}-k-1,m'_{1p}-2k-1);\]
\item  
if $m'_{1p} = 2k+1, n'_{6p} \ge  k+1, n'_{3p} = 1$ 
  or $n'_{6p}=k, m'_{1p} > 2k, n'_{3p}=1$, then
  \[(n_{1p},n_{3p},n_{6p},m_{1p}) :=
  (n'_{1p}+3k+2,0,n'_{6p}-k,m'_{1p}-2k-1);\] 
\item 
otherwise, with
  $n'_{6p}=k, m'_{1p}>2k, n'_{3p}=0$ or $m'_{1p}=2k, n'_{6p} \ge k$:
  \[(n_{1p},n_{3p},n_{6p},m_{1p}) :=
  (n'_{1p}+3k,n'_{3p},n'_{6p}-k,m'_{1p}-2k).\]
\end{itemize}
Conversely, define 
\[
\begin{array}{cccc}\Phi : & \TT_1        &  \to & \TT_2\\ 
                          &(\xi_i,\tau_i)& \mapsto &  (\xi'_i,\tau'_j),
\end{array}
\] where $n'_{ip} := n_{ip}$ for $i \in \{2,\ell,4,5\}$ and
  $m'_{jp} := m_{jp}$ for $j \in \{2,\ell\}$, and the values of
  $n'_{1p},n'_{3p},n'_{6p},m'_{1p}$ depend on $n_{1p}$ modulo $3$ and
  whether $n_{3p}$ is $0$ or $1$:
\begin{itemize}  
\item  
If $n_{1p} \in \{3k+1, 3k+2\}$ and $n_{3p} = 1$, then 
  \[(n'_{1p},n'_{3p},n'_{6p},m'_{1p}) :=
  (n_{1p}-3k-1,0,n_{6p}+k+1,m_{1p}+2k+1);\]
\item  
if $n_{1p} = 3k+2$ and $n_{3p} = 0$, then
  \[(n'_{1p},n'_{3p},n'_{6p},m'_{1p}) :=
  (n_{1p}-3k-2,1,n_{6p}+k,m_{1p}+2k+1);\]
\item  
otherwise, with $n_{1p} \in \{3k,3k+1\}$:
  \[(n'_{1p},n'_{3p},n'_{6p},m'_{1p}) :=
  (n_{1p}-3k,n_{3p},n_{6p}+k,m_{1p}+2k).\]
\end{itemize}  
  It is not difficult to check that $\Phi$ and $\Phi'$ are
  well-defined, that $(\xi_i, \tau_j) \in \TT(\ZZ)$ and $(\xi'_i,
  \tau'_j) \in \TT(\ZZ)$ correspond to the same point $\xx \in U(\QQ)$
  under the map $\Psi$, and that $\Phi$ and $\Phi'$ are inverse to
  each other.
\end{proof}

Together, Lemma~\ref{lem:bijection_different_coprime} and
Lemma~\ref{lem:bijection_different_correct} prove
Proposition~\ref{prop:bijection}.

\section{Congruences}\label{sec:congruences}

We will use the following results from Chapter 3 of \cite{math.NT/0412086} on
the number of solutions of linear and quadratic congruences.

Let $\eta(a;q)$ be the number of positive integers $n \le q$ such that
$\congr{n^2} a q$. Then by equation (3.1) of \cite{math.NT/0412086}, we have
for any $q \in \ZZp$
\begin{equation}\label{eq:bound_eta}
  \eta(a;q) \le 2^{\omega(q)},
\end{equation}
where $\omega(q)$ is the number of distinct prime factors of $q$.
Let $\vartheta$ be an arithmetic function such that 
\[\sum_{d=1}^\infty |(\vartheta \ast \mu)(d)| < \infty,\]
where $\vartheta \ast \mu$ is the usual Dirichlet convolution.

\begin{lemma}\label{lem:sum_arithmetic}
Let $a,q \in \ZZ$ such that $q >
  0$ and $\hcf(a,q)=1$. Then
  \[\sum_{\substack{n \le t\\\congr n a q}} \vartheta(n) = \frac t q
  \sum_{\substack{d = 1\\\hcf(d,q)=1}}^\infty (\vartheta \ast
    \mu)(d) + O\left(\sum_{d=1}^\infty |(\vartheta \ast
      \mu)(d)|\right).\]
\end{lemma}

\begin{proof}
  This is the case $\kappa=0$ of Lemma 2 of \cite{math.NT/0412086}.
\end{proof}

Let $\psi(t) = \{t\} - 1/2$ where $\{t\}$ is the fractional part of
$t \in \RR$. Let $\tpsi(t) = \psi(t)+1$ for $t \in \ZZ$ and $\tpsi(t)
= \psi(t)$ otherwise.

\begin{lemma}\label{lem:sum_congruence}
  Let $a,q \in \ZZ$, where $q > 0$ and $\hcf(a,q)=1$. Let $b_1,b_2 \in \RR$
  with $b_1 \le b_2$. Then \[\#\{n \mid b_1 \le n \le b_2, \congr n a q\} =
  \frac {b_2-b_1} q + r(b_1,b_2;a,q),\] where
  \[r(b_1,b_2;a,q) = \tpsi\left(\frac{b_1-a}{q}\right) -
  \psi\left(\frac{b_2-a}{q}\right).\]
\end{lemma}

\begin{proof}
  This is a slight generalization of Lemma 3 of \cite{math.NT/0412086}.
\end{proof}

\begin{lemma}\label{lem:sum_square}
  Let $\ep > 0$ and $t \ge 0$. Let $a,q \in \ZZ$ such that $q > 0$ and
  $\hcf(a,q)=1$. Then
  \[\sum_{\substack{0 \le \varrho < q \\ \hcf(\varrho,q)=1}} \psi\left(\frac{t
      - a\varrho^2}{q}\right) \lle q^{1/2+\ep} \qquad\text{and}\qquad
  \sum_{\substack{0 \le \varrho < q \\ \hcf(\varrho,q)=1}} \tpsi\left(\frac{t
      - a\varrho^2}{q}\right) \lle q^{1/2+\ep}.\]
\end{lemma}

\begin{proof}
  For $\psi$, this is Lemma 5 of \cite{math.NT/0412086}.
  
  Note that if $\congr{t}{a\varrho_i^2}{q}$ for $i \in \{1,2\}$, then
  $\congr{\varrho_1}{\pm \varrho_2}{q}$, which implies that there are at most
  two different values for $\varrho$ with $0 \le \varrho < q$ such that
  $(t-a\varrho^2)/q$ is integral. Therefore, the sum for $\tpsi$ differs from
  the one for $\psi$ by at most $2$.
\end{proof}

\section{Summations}\label{sec:summations}

Note that $\tau_\ell$ is determined uniquely by $T$ and the other
variables once a certain congruence is fulfilled. Therefore, our
strategy is first to compute the number of possible $\ttwo$ depending
on $\tone, \xi_i$ such that there exists a unique $\tell$ satisfying
$T$. By summing over $\tone$, the number of possible $\tau_i$ is then
computed depending on $\xi_i$. The summation over the variables
$\xi_i$ is finally handled using the height zeta function.

Let \[X_1 = (B \base {-2} 0 {-1} 3 2 1 0)^{1/3}, \qquad X_2 = (B \base
0 {-1} 0 3 2 1 0)^{1/2},\]and \[X_0 = (B^{-1} \base 2 3 4 3 4 5
6)^{1/6}.\] Then the height conditions $|x_i| \le B$ lift to
\begin{equation}\label{eq:height_ttwo}
\left|\left(\frac{\tau_2}{X_2}\right)X_0^3\right|
\le 1,\qquad \left|\left(\frac{\tau_2}{X_2}\right)^2 +
  \left(\frac{\tau_1}{X_1}\right)^3\right| \le 1,
\end{equation}
and
\begin{equation}\label{eq:height_xi}
  |X_0| \le 1,
\end{equation}
and
\begin{equation}\label{eq:height_tone}
  \left| \left(\frac{\tau_1}{X_1}\right)X_0^4
  \right| \le 1,  
\end{equation}
respectively.

Using Proposition~\ref{prop:bijection}, we can now relate
the counting of rational points of bounded height on
$U \subset \S$ to a count on the torsor.

\begin{lemma}\label{lem:counting_torsor}
  We have
  \[\NU(B) = \#\Bigg\{(\xi_i,\tau_i) \in \TT(\ZZ) \Bigg| 
  \begin{aligned}
    &\text{\eqref{eq:hcf_xi}, \eqref{eq:hcf_tone}, 
      \eqref{eq:hcf_ttwo_tell}, \eqref{eq:height_ttwo}, 
      \eqref{eq:height_xi}, \eqref{eq:height_tone} hold},\\
    &\xi_i > 0
  \end{aligned}
  \Bigg\}.\]
\end{lemma}

\subsection{Summation over  $\ttwo$ and $\tell$}

Let $\xi_i,\tau_1$ satisfy the coprimality conditions
\eqref{eq:hcf_xi}, \eqref{eq:hcf_tone} and the height conditions
\eqref{eq:height_xi} and \eqref{eq:height_tone}. Let $\N =
\N(\xi_i,\tau_1)$ be the number of pairs $(\ttwo,\tell)$ such that
\eqref{eq:hcf_ttwo_tell}, \eqref{eq:height_ttwo} and the torsor
equation $T$ are fulfilled.
Then a M\"obius inversion gives \[\N = \sum_{\kell | \Cell}
\mu(\kell)\N_{\kell},\] where $\N_\kell$ has the same definition as $\N$
except that $\hcfell$ is replaced by the condition $\kell | \tell$,
and $T$ is replaced by 
\[
T_\kell(\xi_i, \tau_i, \kell) =
\kell\Tell+\Ttwo+\Tone = 0.
\]  
Note that $\ttwo$ together with
$T_\kell$ defines $\tell$ uniquely once a certain congruence is
fulfilled. Therefore, 
\[\N_\kell = \#\{\ttwo \mid \hcftwo,
\text{~\eqref{eq:height_ttwo} holds}, \congr
{-\Ttwo}{\Tone}{\kell\Fell} \}.\]
Note that \[\hcf(\Tone,\kell\Fell) 
=\hcf(\Tone,\Cell,\Ttwo) = \hcf(\tone,\xi_6,\Ttwo) = 1\] and
\[\hcf(\Ftwo,\kell\Fell) = \hcf(\Ftwo,\kell) = \hcf(\Ftwo,\Cell,\Tone)
= 1.\]

Therefore, it is enough to sum over all $\kell | \Cell$ with
$\hcf(\kell, \tone\xi_1\xi_2\xi_3) = 1$, and since $\kell|\Cell$
implies $\hcf(\kell,\tone) = 1$, we reduce to $\hcf(\kell,
\xi_1\xi_2\xi_3) = 1$.

This implies that there is a unique integer $\varrho$ satisfying
$\boundrho$ and $\hcfrho$ such that \[\congrrho \text{ and
}\congrrhotwo .\]

We have \[\N = \sum_{\substack{\kell|\Cell\\
    \hcf(\kell,\xi_1\xi_2\xi_3) = 1}}\mu(\kell) \sum_{\substack{\boundrho\\
    \congrrhotwo\\ \hcfrho}} \N_\kell(\varrho)\] where \[\N_\kell(\varrho) =
\#\{\ttwo \mid \hcftwo, \text{~\eqref{eq:height_ttwo} holds},
\congrrho\}.\]
We also know that $\hcf(\varrho\tone\xi_1, \kell\Fell) = 1$.  Now we
can apply Lemma~\ref{lem:sum_arithmetic} to the characteristic function
\[\chi(n) = \begin{cases}
  1, &\text{if $\hcf(n,\Ctwo)=1$}\\
  0, &\text{else.}
\end{cases}\]
Since \[\sum_{\substack{d=1\\ \hcf(d,\kell\Fell)=1}}^\infty
\frac{(\chi \ast \mu)(d)}{d} = \prod_{\substack{p\mid\Ctwo\\p \nmid
    \kell\Fell}} \left(1-\frac 1 p\right) = \prod_{p|\Ctwo}
\left(1-\frac 1 p\right) = \phi^*(\Ctwo),\] where we use
$\hcf(\Ctwo,\kell\Fell) = 1$ and the notation $\phi^*(n) := \phi(n)/n$
as in (5.10) of \cite{math.NT/0412086}, we conclude
\[\N_\kell(\varrho) =
\frac{\phi^*(\Ctwo)X_2}{\kell\Fell} \gone(\tone/X_1, X_0) +
O(2^{\omega(\Ctwo)}),\] where $X_2 \gone(\tone/X_1,X_0)$ gives the total
length of the intervals in which $\ttwo$ must lie by \eqref{eq:height_ttwo},
with 
\begin{equation}\label{eq:integral_gone}
  \gone(u,v) = \int_{\{t \in \RR \mid |tv^3| \le 1, |t^2+u^3| \le 1\}} 1
\dd t.
\end{equation}

By equation \eqref{eq:bound_eta}, the number of integers $\varrho$
with $\boundrho$ such that $\hcfrho$ and $\congrrhotwo$ is at most
\[ \eta(\xi_2\xi_3\tone;\kell\Fell) \le 2^{\omega(\kell\Fell)} \le
2^{\omega(\xi_\ell\xi_4\xi_5\xi_6)}.\] This gives as the first step
towards the proof of Theorem \ref{thm:main}:

\begin{lemma}
\[\N = \frac{X_2}{\Fell}\gone(\tone/X_1, X_0)\Sigma(\xi_i,\tau_1) +
O(2^{\omega(\Ctwo)+\omega(\xi_\ell\xi_4\xi_5\xi_6)+\omega(\Cell)}),\]
where
\[\Sigma(\xi_i,\tau_1) =
\phi^*(\Ctwo)\sum_{\substack{\kell|\Cell\\ \hcf(\kell,\xi_1\xi_2\xi_3) = 1}}
\frac{\mu(\kell)}{\kell} \sum_{\substack{\boundrho\\ \congrrhotwo\\
  \hcf(\varrho,\kell\Fell) = 1}} 1\]
\end{lemma}

Now we show that the error term suffices for Theorem~\ref{thm:main} by summing
it over all the $\xi_i, \tone$ which satisfy the height conditions
\eqref{eq:height_xi} and \eqref{eq:height_tone}; we can ignore the coprimality
conditions \eqref{eq:hcf_xi}, \eqref{eq:hcf_tone}, \eqref{eq:hcf_ttwo_tell}.
We obtain:
\begin{equation*}
  \begin{split}
    &\sum_{\xi_i} \sum_{\tone}
    2^{\omega(\Ctwo)+\omega(\xi_\ell\xi_4\xi_5\xi_6)+\omega(\Cell)}\\
    \ll &\sum_{\xi_i}
    2^{\omega(\Ctwo)+\omega(\xi_\ell\xi_4\xi_5\xi_6)+\omega(\Cell)}
    \frac{X_1}{X_0^4}\\
    = &\sum_{\xi_i}
    2^{\omega(\Ctwo)+\omega(\xi_\ell\xi_4\xi_5\xi_6)+\omega(\Cell)}
    \frac B {\base 2 2 3 1 2 3 4}\\
    \ll &B(\log B)^2 \sum_{\xi_i, i\ne \ell}
    \frac{2^{\omega(\Ctwo)+\omega(\xi_4\xi_5\xi_6)+\omega(\Cell)}}{\base
      2 2 3 0 2 3
      4}\\
    \ll &B(\log B)^2
  \end{split}
\end{equation*}
For $\xell$, we have used the estimate \[\sum_{n \le x} 2^{\omega(n)}
\ll x (\log x)\] together with partial summation.

Therefore, we only need to consider the main term when summing over
$\tau_1, \xi_i$ in order to prove Theorem~\ref{thm:main}.

\subsection{Summation over $\tone$}

For fixed $\xi_2, \dots, \xi_6$ satisfying \eqref{eq:hcf_xi} and
\eqref{eq:height_xi}, we sum over all $\tone$ satisfying the
coprimality condition \eqref{eq:hcf_tone} and the height condition
\eqref{eq:height_tone}. Let \[\N' = \N'(\xi_i) = \frac{X_2}\Fell
\sum_{\substack{\tone, \text{~\eqref{eq:height_tone} holds}\\
    \hcfone}} \gone(\tone/X_1, X_0)\Sigma(\xi_i,\tone).\]

First, we find an asymptotic formula for \[\NN(b_1,b_2) = \phi^*(\Ctwo)
\sum_{\substack{\kell|\Cell\\ \hcf(\kell,\xi_1\xi_2\xi_3) =
    1}}\frac{\mu(\kell)}{\kell}\sum_{\substack{\boundrho\\
    \hcfrho}} \N'_\kell(\varrho,b_1,b_2),\] where
\[\N'_\kell(\varrho,b_1,b_2) = \#\left\{\tone \in [b_1,b_2]
  \biggl\lvert \begin{aligned}&\hcfone,\\ &\congrrhotwo
\end{aligned}
\right\}.\]

We have $\hcf(\varrho^2\xi_2, \kell\Fell) = 1$. Therefore we can
replace the condition $\hcfone$ by $\hcf(\tone,\xi_2\xi_3\xi_6) = 1$
in the definition of $\N'_\kell(\varrho,b_1,b_2)$
and perform another M\"obius inversion to
obtain 
\begin{equation*}
  \begin{split}
    \NN(b_1,b_2) = &\phi^*(\Ctwo)\sum_{\substack{\kell|\Cell\\ 
        \hcf(\kell,\xi_1\xi_2\xi_3) = 1}}
    \frac{\mu(\kell)}\kell \sum_{\substack{k_1|\xi_2\xi_3\xi_6\\
        \hcf(k_1,\kell\Fell)=1}} \mu(k_1) \sum_{\substack{\boundrho\\
        \hcfrho}}\\
    &\N'_{\kell,k_1}(\varrho, b_1,b_2),
  \end{split}
\end{equation*}
where
\[\N'_{\kell,k_1}(\varrho,b_1,b_2) = \#\{\tone \in [b_1/k_1, b_2/k_1] 
\mid \congr{-\varrho^2\xi_2}{k_1\tone\xi_3}{\kell\Fell} \}.\] Note
that we must only sum over the $k_1$ with $\hcf(k_1,\kell\Fell)$
because of $\hcf(\varrho^2\xi_2,\kell\Fell)=1$.

Let $a = a(\xi_i,k_1,\kell)$ be the unique integer such that $0 < a
\le \kell\Fell$ and \[\congr{-\xi_2}{k_1 a \xi_3}{\kell\Fell}.\] Then
$\congr{-\varrho^2\xi_2}{k_1 \tone\xi_3}{\kell\Fell}$ if and only if
$\congr{\tone}{a \varrho^2}{\kell\Fell}$. Since
$\hcf(\xi_2,\kell\Fell) = 1$, we have $\hcf(a, \kell\Fell) = 1$. Using
Lemma~\ref{lem:sum_congruence}, we conclude
\[\N'_{\kell,k_1}(\varrho,b_1,b_2) = \frac
{b_2-b_1}{k_1\kell\Fell} + r(b_1/k_1,b_2/k_1,a\varrho^2,\kell\Fell)\]
where, by definition of $r$,
\[r(b_1/k_1,b_2/k_1, a \varrho^2, \kell\Fell) =
\tpsi\left(\frac{b_1/k_1-a\varrho^2}{\kell\Fell}\right) -
\psi\left(\frac{b_2/k_1-a\varrho^2}{\kell\Fell}\right).\]

Let \[\vartheta = \vartheta(\xi_i) =
\begin{cases}
  \phi^*(\Ctwo)\phi^*(\Cone)
  \frac{\phi^*(\Cell)}{\phi^*(\hcf(\Cell,\xi_1\xi_2\xi_3))},
  &\text{\eqref{eq:hcf_xi} holds,}\\
  0, &\text{otherwise.}
\end{cases}
\]
Then for any $b_1<b_2$, we have \[\NN(b_1,b_2) = \vartheta(\xi)\cdot
(b_2-b_1) + \R(b_1,b_2)\] where
\begin{equation*}
  \begin{split}
    \R(b_1,b_2) = &\phi^*(\Ctwo)\sum_{\substack{\kell|\Cell\\\hcfell}} 
    \frac{\mu(\kell)}\kell \sum_{\substack{k_1|\xi_2\xi_3\xi_6\\
        \hcf(k_1,\kell\Fell)=1}} \mu(k_1)
    \sum_{\substack{\boundrho\\\hcfrho}}\\
    &r(b_1/k_1,b_2/k_1, a\varrho^2, \kell\Fell).
  \end{split}
\end{equation*}

By partial summation, we obtain \[\N'(\xi_i) =
\frac{\vartheta(\xi_i)X_1X_2}{\Fell} \gtwo(X_0) + R'(\xi_i)\] with 
\begin{equation}\label{eq:integral_gtwo}
  \gtwo(v) = \int_{\{u \in \RR \mid |uv^4| \le 1\}} \gone(u,v) \dd u
\end{equation}
and
\[R'(\xi_i) = \frac{-X_2}{\Fell}
\int_{-X_0^{-4}}^{X_0^{-4}}
(D_1\gone)(u,X_0)\R(-X_1/X_0^4,X_1 u) \dd u\]
where $D_1\gone$ is the derivation of $\gone$ with respect to the first
variable.

\begin{lemma}\label{lem:summation_tone}
  For any $\xi_i$ as in \eqref{eq:hcf_xi}, \eqref{eq:height_xi}, we have
  \[\N'(\xi_i) = \frac{\vartheta(\xi_i)X_1X_2}{\Fell}\gtwo(X_0) + R'(\xi_i)\]
  where the error term $R'(\xi_i)$ satisfies \[\sum_{\xi_i,
    \text{\eqref{eq:hcf_xi}, \eqref{eq:height_xi} holds}} R'(\xi_i) =
  O(B \log B).\]
\end{lemma}

\begin{proof}
  By Lemma~\ref{lem:sum_square}, we have
\begin{equation*}
  \begin{split}
    \R(b_1,b_2) &\lle \phi^*(\Ctwo) \sum_{\substack{\kell|\Cell\\\hcfell}}
    \frac{|\mu(\kell)|}\kell \sum_{\substack{k_1|\xi_2\xi_3\xi_6\\
        \hcf(k_1,\kell\Fell)=1}} |\mu(k_1)| (\kell\Fell)^{1/2+\ep}\\
    &\le \sum_{\substack{\kell|\Cell\\\hcfell}} |\mu(\kell)|
    \sum_{\substack{k_1|\xi_2\xi_3\xi_6\\
        \hcf(k_1,\kell\Fell)=1}} |\mu(k_1)| (\Fell)^{1/2+\ep}\\
    &\le 2^{\omega(\xi_4\xi_5\xi_6)+\omega(\xi_2\xi_3\xi_6)}
    (\Fell)^{1/2+\ep}
  \end{split}
\end{equation*}
Therefore,
\begin{equation*}
  \begin{split}
    R'(\xi_i) &\lle \frac{X_2}{\Fell}
    \int_{-X_0^{-4}}^{X_0^{-4}}
    (D_1\gone)(u,X_0)2^{\omega(\xi_4\xi_5\xi_6)+\omega(\xi_2\xi_3\xi_6)}
    (\Fell)^{1/2+\ep} \dd u\\
    &\ll \frac{X_2}{(\Fell)^{1/2-\ep}}
    2^{\omega(\xi_4\xi_5\xi_6)+\omega(\xi_2\xi_3\xi_6)}
  \end{split}
\end{equation*}
Summing this over all $\xi_i \le B$, we get using \eqref{eq:height_xi}
\begin{equation*}
  \begin{split}
    \sum_{\substack{\xi_i \le B\\\text{\eqref{eq:hcf_xi}, \eqref{eq:height_xi}
        hold}}} R'(\xi_i) &\lle \sum_{\xi_i \le B}
    \frac{X_2}{X_0^3(\Fell)^{1/2-\ep}}
    2^{\omega(\xi_4\xi_5\xi_6)+\omega(\xi_2\xi_3\xi_6)}\\
    &=\sum_{\xi_i \le B} \frac B {\base 1 2 2 {3/2-3\ep} {2-2\ep}
      {5/2-\ep} 3}
    2^{\omega(\xi_4\xi_5\xi_6)+\omega(\xi_2\xi_3\xi_6)}\\
    &\ll \sum_{\xi_i, i\ne 1} \frac {B\log B}{\base 0 2 2
      {3/2-3\ep} {2-2\ep} {5/2-\ep}
      3}2^{\omega(\xi_4\xi_5\xi_6)+\omega(\xi_2\xi_3\xi_6)}\\
    &\ll B\log B\qedhere
  \end{split}
\end{equation*}
\end{proof}

\subsection{Summation over $\xi_i$}

Define \[\Delta(n) = B^{-5/6}\sum_{\xi_i, \base 2 3 4 3 4 5 6 = n}
\frac{\vartheta(\xi_i) X_1X_2}\Fell.\]

We sum $\N'(\xi_i)$ in Lemma~\ref{lem:summation_tone} over the seven
variables $\xi_i$ such that the coprimality conditions
\eqref{eq:hcf_xi} and the height condition \eqref{eq:height_xi} hold.
Note that the definition of $\vartheta(\xi_i)$ ensures that the main
term of $\N'(\xi_i)$ is zero if \eqref{eq:hcf_xi} is not satisfied. In
view of Lemma~\ref{lem:counting_torsor}, this implies:

\begin{lemma}\label{lem:NUB}
  We have \[\NU(B) = B^{5/6} \sum_{n \le B} \Delta(n) \gtwo((n/B)^{1/6}) +
  O(B(\log B)^2).\] 
\end{lemma}

\section{Proof of Manin's Conjecture for $\S$}\label{sec:main_thm_proof}

Our argument is similar to \cite{math.NT/0412086}.
We need to estimate \[M(t) := \sum_{n \le t}
\Delta(t)\] for $t > 1$. Therefore, we consider the Dirichlet series
$F(s) := \sum_{n=1}^\infty \Delta(n)n^{-s}$.

Observing \[\frac{X_1X_2}{\Fell} = \frac{B^{5/6}(\base 2 3 4 3 4 5
  6)^{1/6}}{\base 1 1 1 1 1 1 1},\] we get \[F(s+1/6) = \sum_{\xi_i}
\frac{\vartheta(\xi_i)}{\xi_1^{2s+1}\xi_2^{3s+1}\xi_3^{4s+1}\xi_\ell^{3s+1}
  \xi_4^{4s+1}\xi_5^{5s+1}\xi_6^{6s+1}},\] and writing $F(s+1/6) =
\prod_p F_p(s+1/6)$ as a product of local factors, we obtain:
\begin{equation*}
\begin{split}
  F_p(s+1/6) = &1 + \frac{(1-1/p)^2}{(p^{\lambda_6s+1}-1)}\bigg(
  \frac{p^{\lambda_1s+1}}{p^{\lambda_1s+1}-1} +
  \frac{p^{\lambda_1s+1}p^{\lambda_6s+1}}{p^{\lambda_3s+1}
    (p^{\lambda_1s+1}-1)}\\
  &+ \frac{p^{\lambda_6s+1}}{(1-1/p)p^{\lambda_2s+1}} +
  \frac{1}{p^{\lambda_\ell s+1}-1}
  + \frac{p^{\lambda_\ell s+1}p^{\lambda_6s+1}}
  {p^{\lambda_4s+1}(p^{\lambda_\ell s+1}-1)}\\ &+ \frac{p^{\lambda_\ell
      s+1} p^{\lambda_6s+1}}{p^{\lambda_5s+1}(p^{\lambda_\ell
      s+1}-1)}\bigg) + \frac{1-1/p}{p^{\lambda_1s+1}-1} +
  \frac{1-1/p}{p^{\lambda_\ell s+1}-1}
\end{split}
\end{equation*}
for any prime $p$.

Since $1/p^{\lambda_i s +1}=O_\ep(1/p^{1/2+\ep})$ for $s \in \HH := \{s \in
\CC \mid \re(s) \ge -1/12 + \ep\}$ and $i \in \{1,2,3,\ell,4,5,6\}$, we have
\[F_p(s+1/6) = 1 + \sum_i \frac{1}{p^{\lambda_i s + 1}} +
O_\ep\left(\frac{1}{p^{1+\ep}}\right)\] for $s \in \HH$, and defining
\[E(s):=\prod_i\zeta(\lambda_is+1) =
  \zeta(2s+1)\zeta(3s+1)^2\zeta(4s+1)^2\zeta(5s+1)\zeta(6s+1),\]we have
\[\frac 1{E_p(s)} = 1 - \sum_i \frac 1{p^{\lambda_i s +
    1}} + O_\ep\left(\frac{1}{p^{1+\ep}}\right)\] for $s \in \HH$.
Define \[G(s) := F(s+1/6)/E(s)\] for $\re(s)>0$. Then $G$ has a holomorphic
and bounded continuation to $\HH$. Note that
\[G(0) = \prod_p \left(1-\frac 1 p\right)^7 \left(1+\frac 7 p + \frac 1
  {p^2}\right),\] and that for $s \to 0$, we have
\[E(s) = \frac 1{\prod_i \lambda_i}
  s^{-7} + O(s^{-6} ).\]

Therefore, the residue of $F(s)t^s/s$ at $s = 1/6$ is \[\Res(t) =
\frac{6 G(0)t^{1/6}Q_1(\log t)}{6! \cdot \prod_i \lambda_i}\] for some
monic polynomial $Q_1$ of degree $6$.

\begin{lemma}
  $M(t) = \omega_0 \alpha(\tS)\cdot 6 t^{1/6}Q_1(\log t) +
  O_\ep(t^{1/6-1/24+\ep})$.
\end{lemma}

\begin{proof}
  Integrating Perron's formula for $M(t)$ over $t$, we have
  \[\int_0^t M(u) \dd u = \frac 1{2\pi i}
  \int_{1/6+\ep-i\infty}^{1/6+\ep+i\infty} F(s) \frac
  {t^{s+1}}{s(s+1)}\dd s\] for $t>1$ and $\ep > 0$.
  
  We apply Cauchy's residue theorem to the rectangle with vertices
  \[1/12+\ep-iT, 1/12+\ep+iT, 1/6+\ep+iT, 1/6+\ep-iT,\]
  for some $T>1$, where $\ep > 0$ is sufficiently small.
  
  By the convexity bound \[\zeta(1+\sigma+i\tau) \lle
  |\tau|^{-\sigma/3+\ep}\] for any $\sigma\in [-1/2,0)$, we have
  \begin{equation}\label{eq:convexity}
    F(1/6+\sigma+i\tau) \ll E(\sigma+i\tau) \lle |\tau|^{-9\sigma+\ep}
  \end{equation}
  for any $\sigma \in [-1/12+\ep,0)$, using $(\sum_i \lambda_i)/3 = 9$
  and that $G(\sigma+i\tau)$ is bounded.

  For the ray going down from $1/6+\ep-iT$, we get
  \begin{equation*}
    \begin{split}
      \left|\int_{1/6+\ep-i\infty}^{1/6+\ep-iT} \frac{F(s)t^{s+1}}{s(s+1)}
        \dd s\right| &\le \int_{-\infty}^T
      \frac{|F(1/6+\ep+i\sigma)||t^{7/6+\ep+i\sigma}|}
      {|(1/6+\ep+i\sigma)(7/6+\ep+i\sigma)} \dd \sigma\\ &\ll
      t^{7/6+\ep} \int_{-\infty}^T \frac 1{|\sigma|^2} \dd
      \sigma
      \\ &\ll t^{7/6+\ep}T^{-1}
    \end{split}
  \end{equation*}
  where we use that $F(s)$ is bounded for $\re(s) \ge 1/6+\ep$.
  Integrating from $1/6+\ep+iT$ to $1/6+\ep+i\infty$ gives
  the same result.

  For the lower edge, we estimate
  \begin{equation*}
    \begin{split}
      \left|\int_{1/12+\ep-iT}^{1/6+\ep-iT} \frac{F(s)t^{s+1}}{s(s+1)}
        \dd s\right| &\le \int_{-1/12+\ep}^\ep
        \frac{|F(1/6+\sigma-iT)||t^{7/6+\sigma-iT}|}
        {|(1/6+\sigma-iT)(7/6+\sigma-iT)|} \dd \sigma\\
        &\lle \frac{T^{9/12+\ep} t^{7/6+\ep}}{T^2},
    \end{split}
  \end{equation*}
  because \eqref{eq:convexity} gives a bound for $-1/12+\ep \le
  \sigma \le -\ep$, $F(s)$ being continuous gives a bound in an
  $\ep$-neighborhood of $1/6-iT$, and the length of the
  integration interval is $1/12$. For the upper
  edge, we obtain the same bound.

  For the edge on the left, we have
  \begin{equation*}
    \begin{split}
      \left|\int_{1/12+\ep-iT}^{1/12+\ep+iT} \frac{F(s)t^{s+1}}{s(s+1)}
        \dd s\right| &\le \int_{-T}^T
      \frac{|F(1/12+\ep+i\sigma)||t^{13/12+\ep+i\sigma}|}
      {|(1/12+\ep+i\sigma)(13/12+\ep+i\sigma)} \dd \sigma\\
      &\lle \int_{-T}^T \frac{|\sigma|^{9/12+\ep}
        t^{13/12+\ep}}{(1+|\sigma|)^2} \dd \sigma\\
      &\ll t^{13/12+\ep}
    \end{split}
  \end{equation*}
  since the integral over $\sigma$ is bounded independently of $T$,
  and using \eqref{eq:convexity} again.
  
  Taking $T=t$, we have proved \[\int_0^t M(u)\dd u = \int_0^t \Res(u)
  \dd u + O_\ep(t^{13/12+\ep}).\] But now \[\frac 1 H \int_{t-H}^t
  M(u) \dd u \le M(t) \le \frac 1 H \int_t^{t+H} M(u) \dd u,\] and for
  $H \le t/3$, both integrals are equal to \[\Res(t) + O_\ep(H t^{5/6}
  (\log t)^6 + H^{-1} t^{13/12+\ep}).\] The proof of the Lemma is
  completed by choosing $H = t^{23/24}$ and noting that $\omega_0 =
  G(0)$ and $\alpha(\tS) = (6! \prod_i \lambda_i)^{-1}$ by the
  definitions of $\omega_0$ and $\alpha(\tS)$ in
  Lemma~\ref{lem:manin_conjecture}.
\end{proof}

By partial summation we conclude 
\begin{equation*}
  \begin{split}
    &\sum_{n\le B} \Delta(n)\gtwo((n/B)^{1/6})\\
    &= \omega_0 \alpha(\tS)\cdot 6 \int_0^B \gtwo(u^{1/6}/B^{1/6}) 
    \frac \dd{\dd u}(u^{1/6}Q_1(\log u)) \dd u + O_\ep(B^{1/6-1/24+\ep})\\
    &= B^{1/6} \omega_0 \alpha(\tS)\cdot 6 \int_0^1\gtwo(v) Q_2(\log B+ 6\log
    v) \dd v + O_\ep(B^{1/6-1/24+\ep})
  \end{split}
\end{equation*}
for some monic polynomial $Q_2$ of degree $6$. Considering definitions
\eqref{eq:integral_gone} and \eqref{eq:integral_gtwo}, note that
 \[\omega_\infty = 6 \int_{\{v \in \RR \mid 0 \le v
   \le 1\}} \gtwo(v) \dd v.\] Together with Lemma~\ref{lem:NUB}, this
 completes the proof of Theorem~\ref{thm:main}.

\bibliographystyle{alpha}

\bibliography{e6}

\end{document}